# Some Nasty Reflexive Groups

Rüdiger Göbel and Saharon Shelah

**Abstract**

In *Almost Free Modules, Set-theoretic Methods*, Eklof and Mekler [5, p. 455, Problem 12] raised the question about the existence of dual abelian groups $G$ which are not isomorphic to $\mathbb{Z} \oplus G$. Recall that $G$ is a dual group if $G \cong D^*$ for some group $D$ with $D^* = \mathrm{Hom}\,(D, \mathbb{Z})$. The existence of such groups is not obvious because dual groups are subgroups of cartesian products $\mathbb{Z}^D$ and therefore have very many homomorphisms into $\mathbb{Z}$. If $\pi$ is such a homomorphism arising from a projection of the cartesian product, then $D^* \cong \ker \pi \oplus \mathbb{Z}$. In all 'classical cases' of groups $D$ of infinite rank it turns out that $D^* \cong \ker \pi$. Is this always the case? Also note that reflexive groups $G$ in the sense of H. Bass are dual groups because by definition the evaluation map $\sigma : G \longrightarrow G^{**}$ is an isomorphism, hence $G$ is the dual of $G^*$. Assuming the diamond axiom for $\aleph_1$ ($\diamondsuit_{\aleph_1}$) we will construct a reflexive torsion-free abelian group of cardinality $\aleph_1$ which is not isomorphic to $\mathbb{Z} \oplus G$. The result is formulated for modules over countable principal ideal domains which are not field.

## 1 Introduction

Let $R$ be a countable principal ideal domain with $1 \neq 0$ and not a field. If $S = R \setminus \{0\}$ is enumerated by $s_n$ $(n \in \omega)$ such that $s_0 = 1$, then we let $q_n = \prod_{i<n} s_i$. The $q_n$'s constitute a divisibility chain with $q_{n+1} = q_n s_n$ and $q_0 = 1$. The principal ideals $q_n R$ define a neighborhood basis of $0 \in R$ for the $R$-topology of $R$ which is Hausdorff because

---

This work is supported by the project No. G-545-173.06/97 of the German-Israeli Foundation for Scientific Research & Development
secondary: 03E05, 03E35
GbSh 568 in Shelah's list of publications

$\bigcap_{n \in \omega} q_n R = 0$. The $R$-completion of $R$ is denoted by $\widehat{R}$ and obviously $|\widehat{R}| = 2^{\aleph_0}$, see also Göbel, May [7]. If $F$ is a free $R$-module, then similarly $F \subseteq \widehat{F}$ and $F$ is pure and dense in $\widehat{F}$. Recall that $F \subseteq_* \widehat{F}$ is pure if and only if $\widehat{F} q_n \cap F \subseteq F q_n$ for all $n \in \omega$. Also $F \subseteq \widehat{F}$ is dense if and only if $\widehat{F}/F$ is divisible. The fact that $|\widehat{R}| = 2^{\aleph_0}$ is reflected in an easy

**Observation 1.1** *If $0 \neq r_n \in R$ for all $n \in \omega$ then we can find $p_n \in \{0, r_n\}$ such that*

$$\sum_{n \in \omega} p_n q_n \in \widehat{R} \setminus R.$$

See again [7].

We will use topological arguments and the prediction principle $\diamondsuit_{\aleph_1}$ which holds in many models of set theory, in particular in Gödel's universe $V = L$, to answer a problem in the book by Eklof and Mekler [5, p. 455, Problem 12] concerning dual modules. Recall that a module $G$ is a dual module if $G \cong D^*$ where $D^* = \text{Hom}(D, R)$ is the dual of the $R$-module $D$. There is a large range of recent literature on dual modules which can be looked up in [5]. Strongly related with dual modules are the well-known evaluation maps

$$\sigma = \sigma_D : D \longrightarrow D^{**} \quad (d \longrightarrow \sigma(d))$$

where $\sigma(d) \in D^{**}$ is the homomorphism defined by evaluation

$$\sigma(d) : D^* \longrightarrow R \quad (\varphi \longrightarrow \varphi(d)).$$

Investigating questions on homology Bass [1, p. 476] introduced the notion *reflexivity* and *torsion-less*. Now $D$ is *reflexive* if and only if the evaluation map $\sigma_D$ is an isomorphism. ($D$ is torsion-less if this map is an injection.) Trivially reflexive modules $D$ are dual modules as $D \cong (D^*)^*$. Moreover it is obvious that dual modules $G \cong D^*$ are ( isomorphic to ) submodules of the cartesian product $R^D$. An $R$-module is $\aleph_1$-free if all its countable submodules are free and recall from a result of Specker (see Fuchs [6]) that cartesian products $R^\kappa$ are $\aleph_1$-free. Hence dual modules are $\aleph_1$-free as well. In particular reflexive modules $G \cong (G^*)^*$ are dual modules hence $\aleph_1$-free and the example $G$ we want to construct must be $\aleph_1$-free. This is also a warning that anticipated results may depend on the set theory in use. We have seen that reflexive modules are dual modules, hence the following theorem – the main target of this paper – provides a strong negative answer to the problem in [5].

**Theorem 1.2** *(ZFC $+ \diamondsuit_{\aleph_1}$) If $R$ countable principal ideal domain which is not a field, then there is a reflexive $R$-module $G$ of cardinality $\aleph_1$ such that $G \not\cong R \oplus G$.*

Using a remark on Γ-invariants from Section 3 we can construct a family of pairwise non-isomorphic examples which has size $2^{\aleph_1}$.

It is also remarkable that Eda's interesting examples of dual groups arising from a theory of continuous functions do not share the main property of the theorem, see Eda [3] and Eda, Otah [4].

This paper is also basic for proving the same theorem in ZFC assuming the weaker special continuum hypothesis CH only. The present work, so to speak, represents the 'local case' for [8]. The question whether the result holds in any model of ZFC remains open. In another paper [9] however we will show that $(\diamondsuit_{\aleph_1})$ or CH need not hold for showing Theorem 1.2. We are also able to derive the statement assuming Martin's Axiom (and e.g. the negation of CH).

We will work in the category of free modules with bilinear forms leading to some torsion-free module with bilinear form, which resembles the Hahn-Banach-Theorem from functional analysis.

## 2   Bilinear Forms on Free $R$-modules

Let $R$ be the principal ideal domain discussed in Section 1. A *bilinear form*

$$\Phi : G \oplus H \longrightarrow R$$

is a map with domain $\mathrm{Dom}\, \Phi = G \oplus H$ an $R$-module such the following two conditions hold. If $g \in G$, then

$$\Phi(g,\ ) : H \longrightarrow R$$

is an $R$-homomorphism (we say a homomorphism for short) and dually if $h \in H$, then

$$\Phi(\ ,h) : G \longrightarrow R$$

is a homomorphism as well.

In this definition we call $G$ the left and $H$ the right part of $\Phi$. If $H^* = \mathrm{Hom}\,(H,R)$ is the dual module, then $\Phi(g,\ ) \in H^*$, $\Phi(\ ,h) \in G^*$ and

$$\Phi(G,\ ) \subseteq H^*,\ \Phi(\ ,H) \subseteq G^*.$$

We will also consider a particular class $\mathfrak{F}$ of such bilinear forms.

**Definition 2.1** *Say that the bilinear form $\Phi$ belongs to $\mathfrak{F}$ if $\Phi : G \oplus H \longrightarrow R$ satisfies the following conditions.*

*(i) $G$ and $H$ are free $R$-modules of countable rank.*

(ii) $\Phi$ *is non-degenerative, that is* $\Phi(g,\ )=0$ *only if* $g=0$ *and dually* $\Phi(\ ,h)=0$ *only if* $h=0$ *for any* $g\in G$ *and* $h\in H$.

(iii) $\Phi$ *preserves purity, that is* $\Phi(g,\ )\in_* H^*$ *if* $g\in_* G$ *and dually* $\Phi(\ ,h)\in_* G^*$ *if* $h\in_* H$.

For brevity we will call the conditions in Definition 2.1 the membership conditions for $\Phi$ (to belong to $\mathfrak{F}$). Also recall that $g\in_* G$ denotes a pure element or equivalently $gR$ is a pure submodule of $G$.

If $G=H=\bigoplus_{n\in\omega}e_nR$ and $g=\sum_{n\in\omega}e_ng_n\in G$ is the usual direct sum representation of $g$ with $g_n\in R$ and $g_n=0$ for almost all $n\in\omega$, then the natural scalar product

$$\Phi(g,h)=\sum_{n\in\omega}g_nh_n\in R$$

is a bilinear form and it is easy to see that $\Phi\in\mathfrak{F}$, hence $\mathfrak{F}$ is not empty. We define an ordering on $\mathfrak{F}$ by taking $\Phi\subseteq\Phi'$ in $\mathfrak{F}$ for two bilinear forms if $\Phi'$ extends $\Phi$ and $\mathrm{Dom}\,\Phi$ is a pure submodule of $\mathrm{Dom}\,\Phi'$. The proof of the following lemma is easy checking of the definitions above.

**Lemma 2.2** *Let* $\Phi\in\mathfrak{F}$, $\mathrm{Dom}\,\Phi=G\oplus H$ *and* $0\neq f\in_* G^*$, $0\neq g\in_* H^*$ *respectively.*

(i) *There exists* $\Phi'\in\mathfrak{F}$ *such that* $\Phi\subseteq\Phi'$, $\mathrm{Dom}\,\Phi'=G\oplus(H\oplus eR)$ *and* $\Phi'(\ ,e)=f$.

(ii) *There exists* $\Phi'\in\mathfrak{F}$ *such that* $\Phi\subseteq\Phi'$, $\mathrm{Dom}\,\Phi'=(G\oplus eR)\oplus H$ *and* $\Phi'(e,\ )=g$.

The following definitions are crucial for proving the main result.

**Definition 2.3** *If* $\Phi\in\mathfrak{F}$ *with* $\mathrm{Dom}\,\Phi=G\oplus H$ *then* $\varphi\in G^*$ *is* essential *for* $\Phi$ *if for any finite rank summand* $L$ *of* $G$ *and any finite subset* $E$ *of* $H$ *there is* $g\in G\setminus L$ *with*

$$g\varphi\neq 0=\Phi(g,e)\ \text{for all}\ e\in E.$$

The notion *essential* for $\varphi\in H^*$ is dual. Obviously we may assume that $g$ in the Definition 2.3 is such that $L\oplus gR$ is a summand of $G$.

If $g\varphi=0=\Phi(g,e)$ for all $e\in E$ and some finite $E\subseteq H$ then $\varphi$ is a linear combination of the linear combination of the $\Phi(\ ,e)$'s by induction on $|E|$. Hence $\varphi\in G^*$ is essential for $\Phi$ is equivalent to say that $\varphi$ is not in $\Phi(\ ,H)$ *modulo* finite rank in $G$.

**Definition 2.4** *We say that $\Phi \in \mathfrak{F}$ with $\mathrm{Dom}\,\Phi = G \oplus H$ is finitely covered on the left if for any $g \in L^*$, $L \subseteq_* G$ of finite rank we find $h \in H$ with $\Phi(\ ,h) \restriction L = g$. The definition finitely covered on the right is dual and $\Phi$ is finitely covered if it is both finitely covered on the right and on the left.*

We have an immediate

**Corollary 2.5** *If $\Phi \in \mathfrak{F}$ there is $\Phi \subseteq \Phi' \in \mathfrak{F}$ such that $\Phi'$ is finitely covered.*

**Proof** If $L \subseteq_* G$ has finite rank then $L$ is a summand of $G$ and any element of $L^*$ extends to an element of $G^*$. If $f \in_* L^*$ is not of the form $\Phi(\ ,e) \restriction L$ for some $e \in_* H$, then apply Lemma 2.2 to find an extension $\Phi'$ taking care of an extension of $f$ to $G^*$. After countably many steps - taking unions - we find an extension $\Phi \subseteq \Phi_1$ which is finitely covered on the left. Similarly $\Phi_1 \subseteq \Phi_2$ is finitely covered on the right. We proceed this way to find $\Phi_\omega = \bigcup_{n \in \omega} \Phi_n$ which is finitely covered such that $\Phi \subseteq \Phi_\omega \in \mathfrak{F}$.

In the proof of the Corollary 2.5 we used twice the following easy observation.

**Lemma 2.6** *If $\Phi_\alpha$ $(\alpha \in \delta)$ is an ascending, continuous chain of bilinear forms in $\mathfrak{F}$ and $\delta < \aleph_1$, then $\Phi = \cup_{\alpha \in \delta} \Phi_\alpha \in \mathfrak{F}$ and $\Phi_\alpha \subseteq \Phi$ for all $\alpha < \delta$.*

**Proof.** If $\alpha \in \delta$, then $\Phi_\alpha \subseteq \Phi$ as maps and $\mathrm{Dom}\,\Phi_\alpha = G_\alpha \oplus H_\alpha$ is a pure submodule of $\mathrm{Dom}\,\Phi = G \oplus H$ with $G = \cup_{\alpha \in \delta} G_\alpha$, $H = \cup_{\alpha \in \delta} H_\alpha$, hence $\Phi_\alpha \subseteq \Phi$ also as members in $\mathfrak{F}$. Moreover $G \oplus H$ is countable and any pure finite rank submodule belongs to some $G_\alpha \oplus H_\alpha$ by purity, hence it is free. If follows from Pontryagin's theorem that $G \oplus H$ is free, see Fuchs [6, p.93]. The other membership conditions of $\mathfrak{F}$ are automatic, hence $\Phi \in \mathfrak{F}$ as desired.

**First Killing-Lemma 2.7** *Suppose that $\varphi \in G^*$ is essential for $\Phi \in \mathfrak{F}$ with $\mathrm{Dom}\,\Phi = G \oplus H$. Then we find $\Phi \subseteq \Phi' \in \mathfrak{F}$ with*

$$\mathrm{Dom}\,\Phi' = G' \oplus H \text{ and } G' = \langle G, y_0 \rangle_* \subseteq \widehat{G}$$

*for some $y_0 \in \widehat{G}$ such that $\varphi$ does not extend to $G' \longrightarrow R$.*

Remark. By symmetry a similar lemma holds for $\varphi \in H^*$.

**Proof** Let $H = \bigoplus_{i \in \omega} h_i R$ and $G = \bigoplus_{n \in \omega} g_n R$. Inductively we construct elements $g'_n \in G \setminus G_n$ with $G_n = \langle g_i, g'_i : i < n \rangle_* \subseteq G$ such that the following holds:

(i) $\Phi(g'_n, h_i) = 0$ for all $i < n$.

(ii) $\bigoplus_{i<n} g'_i R$ is a direct summand of $G$ - and also $\bigoplus_{i\in\omega} g'_i R$ is a summand of $G$.

(iii) $g'_n \varphi \neq 0$ for all $n \in \omega$.

Suppose $g'_i$ for $i < n$ are constructed accordingly. We must find $g'_n \in G \setminus G_n$ with $(i), (ii), (iii)$. If $h^i = \Phi(\ , h_i)$ for $i \leq n$, then recall that $\varphi$ is essential for $\Phi$ and there is $g'_n \in_* G \setminus G_n$ such that
$$g'_n \varphi \neq 0 = \Phi(g'_n, h_i) = g'_n h^i \text{ for } i \leq n.$$
Hence $(i), (ii), (iii)$ hold and the second part of $(ii)$ follows by an easy support argument.

By Observation 1.1 we now can choose $k_n \in \{0, q_n\}$ such that
$$\sum_{n\in\omega} k_n(g'_n \varphi) \in \widehat{R} \setminus R.$$

Let $y_0 = \sum_{n\in\omega} k_n g'_n \in \widehat{G}$ and
$$G' = \langle G, y_0 \rangle_* \subseteq \widehat{G}.$$
Note that $G' = \langle G, y_s R : s \in \omega \rangle$ with
$$y_s = \sum_{s \leq n \in \omega} k_n (k_s)^{-1} g'_n \text{ if } k_s \neq 0 \text{ and } y_s = 0 \text{ if } k_s = 0.$$
Replacing some of the generators $g_n$ ($n \in \omega$) of $G$ we may assume that the $g'_n$'s and some of the $g_n$'s generate $G$ freely. Then the $y_n \in G'$ can be used to generate some of the $g'_n$'s, hence the $y_n$'s together with some of the mentioned generators of $G$ generate all of $G'$. Obviously these elements (by support) are independent, hence also $G'$ is freely generated.

The homomorphism $\varphi : G \longrightarrow R$ by continuity extends uniquely to $\hat{\varphi} : G' \longrightarrow \widehat{R}$ however
$$y_0 \hat{\varphi} = \sum_{n\in\omega} k_n(g'_n)\varphi \in \widehat{R} \setminus R$$
which is no longer a homomorphism into $R$ hence $\varphi$ does not extend.

Let $\Phi' : G' \oplus H \longrightarrow \widehat{R}$ be the unique extension of $\Phi$. In order to see $\Phi \subseteq \Phi'$ in $\mathfrak{F}$ we must have $\operatorname{Im} \Phi' \subseteq R$. However by $(i)$ we have
$$\Phi'(y_0, h_j) = \Phi'(\sum_{n\in\omega} k_n g'_n, h_j) = \sum_{n\in\omega} k_n \Phi(g'_n, h_j) = \sum_{n<j} k_n \Phi(g'_n, h_j) \in R.$$
Note that $k_s y_s = g + y_0$ for some $g \in G$, hence $\Phi'(k_s y_s, h_j) = \Phi'(g + y_0, h_j) \in R$ and $\Phi'(k_s y_s, h_j) = k_s \Phi'(y_s, h_j) \in \widehat{R}$. By purity we have $\Phi'(y_s, h_j) \in R$ as desired. The membership conditions for $\Phi' \in \mathfrak{F}$ are now easily checked.

**Second Killing-Lemma 2.8** *Let $\Phi \in \mathfrak{F}$ be with $\operatorname{Dom} \Phi = G \oplus H$. Suppose that $\eta : G \longrightarrow G$ is a monomorphism with*

$$G = x_0 R \oplus \operatorname{Im} \eta.$$

*Then we can find $\Phi \subseteq \Phi' \in \mathfrak{F}$ such that $\eta$ does not extend to any $\eta'' : G'' \longrightarrow R$ with $\Phi \subseteq \Phi' \subseteq \Phi'' \in \mathfrak{F}$, $\operatorname{Dom} \Phi'' = G'' \oplus H''$ and $G'' = x_0 R \oplus \operatorname{Im} \eta''$.*

**Proof** Inductively we define $x_n \in G$ by

$$G = x_0 R \oplus G\eta = x_0 R \oplus (x_0 R \oplus G\eta)\eta = x_0 R \oplus (x_0\eta) R \oplus G\eta^2 = \cdots = \bigoplus_{i<n} x_i R \oplus G\eta^n.$$

Hence $x_n = x_0 \eta^n$ for all $n \in \omega$. Also let $H = \bigoplus_{n\in\omega} h_n R$ and $G = \bigoplus_{n\in\omega} g_n R$ and let

$g^n \in G^*$ be the homomorphism defined by $g_i g^n = \delta_{i,n}$ $(i \in \omega)$.

For each $n \in \omega$ we want to find $a_i^n \in R$ for $n \leq i \leq 4n$ with the following properties for

$$w = \sum_{i=n}^{4n} x_i a_i^n \text{ and } w' = w\eta = \sum_{i=n}^{4n} x_{i+1} a_i^n.$$

$$w \neq 0, \tag{2.1}$$

$$\Phi(w, h_k) = 0 \text{ and } wg^k = 0 \text{ for all } k < n, \tag{2.2}$$

$$w'g^k = 0 \text{ for all } k < n. \tag{2.3}$$

This is equivalent to say that we seek for a non-trivial solution of a homogeneous system of $3n$ linear equations with $3n+1$ parameters $a_i^n \in R$ $(n \leq i \leq 4n)$. By linear algebra we can find the desired solution $a_i^n \in R$. Similarly we will find a countable family of such elements $w$. Inductively we define an increasing sequence $s_n \in \omega$ $(n \in \omega)$ such that

$$4s_n < s_{n+1} \tag{2.4}$$

and

$$\{x_{s_n}, \ldots, x_{4s_n+1}\} \subseteq \bigoplus_{i<s_{n+1}} g_i R. \tag{2.5}$$

The last two conditions are easily verified. Using (2.4) and (2.1), ... ,(2.3) for $s_n$ in place of $n$ we find an element

$$0 \neq w_n = \sum_{i=s_n}^{4s_n} x_i a_i^n$$

such that

$$\Phi(w_n, h_k) = 0 \text{ and } w_n g^k = 0 \text{ for } k < s_n, \tag{2.6}$$

$$w_n g^k = 0 \text{ for } k < s_n. \tag{2.7}$$

If

$$T_n = \bigoplus_{s_n \leq i < s_{n+1}} g_i R$$

then conditions (2.5), the second part of (2.6) and (2.7) are equivalent to say that

$$w_n, w_n \eta \in T_n \text{ for all } n \in \omega. \tag{2.8}$$

By $[w]$ we denote the *support* of an element $w = \sum x_i a_i$, that is the set $[w] = \{i \in \omega : a_i \neq 0\}$. Let $m_n = \sup [w_n] + 1$ and consider the projection $\pi_n : G \longrightarrow R$ which is defined as follows.

Write $G = \bigoplus_{i<k} x_i R \oplus G\eta^{k+1}$ for some $k \geq m_n$, hence $g \in G$ has a unique (independent of $k$) representation

$$g = \sum_{i \leq k} x_i b_i + g' \text{ with } b_i \in R, g' \in G\eta^{k+1}.$$

Define $g\pi_n = b_{m_n}$ for all $g \in G$. Hence $\pi_n \in G^*$ and clearly

$$w_n \pi_n = 0 \text{ but } w_n \eta \pi_n \neq 0 \text{ for all } n \in \omega \tag{2.9}$$

by the action of $\eta$ on $x_{m_n - 1} \eta = x_{m_n}$. By Observation 1.1 we can find $r_n \in R$ with

$$\sum_{n \in \omega} r_n q_n (w_n \eta) \pi_n \in \widehat{R} \setminus R. \tag{2.10}$$

Define $z = \sum_{n \in \omega} q_n r_n w_n \in \widehat{G}$ and let

$$G' = \langle G, z \rangle_* = \langle G, z_n : n \in \omega \rangle \subseteq \widehat{G}$$

such that $z_0 = z$ and $z_k = \sum_{n \geq k} q_n q_k^{-1} r_n w_n \in \widehat{G}$. By an argument used above we know that $G'$ is a free $R$-module of countable rank. The map $\Phi : G \oplus H \longrightarrow R$ by continuity extends uniquely to $\Phi' : G' \oplus H \longrightarrow \widehat{R}$. We want to see that

$$\Phi' : G' \oplus H \longrightarrow R$$

and must show that $\operatorname{Im} \Phi' \subseteq R$. This will follow from (2.6) and continuity

$$\Phi'(z, h_k) = \Phi'(\sum_{n \in \omega} q_n r_n w_n, h_k) = \sum_{n \in \omega} q_n r_n \Phi(w_n, h_k)$$
$$= \sum_{n \in \omega} q_n r_n \Phi(w_n, h_k) = \sum_{n \leq k} q_n r_n \Phi(w_n, h_k) \in R.$$

Finally we want to extend $\Phi'$ to $\Phi' : G' \oplus H' \longrightarrow R$ with $H' = H \oplus hR$ and put $\Phi'(\ , h) \restriction T_n = \pi_n \restriction T_n$. By $G = \bigoplus_{n \in \omega} T_n$ the map $\Phi'(\ , h) \in G^*$ is well-defined.

By continuity $\Phi' : G \oplus H' \longrightarrow \widehat{R}$ extends uniquely to

$$\Phi' : G' \oplus H' \longrightarrow \widehat{R}$$

and again we must show that $\operatorname{Im} \Phi' \subseteq R$. Note that $\Phi'(w_n, h) = w_n \pi_n = 0$ from above, hence

$$\Phi'(z, h) = \Phi'(\sum_{n \in \omega} q_n r_n w_n, h) = \sum_{n \in \omega} q_n r_n \Phi'(w_n, h) = 0$$

and $\operatorname{Im} \Phi' \subseteq R$.

Also $H'$ is a free $R$-module and $G \oplus H$ is a pure submodule of $G' \oplus H'$, hence $\Phi \subseteq \Phi'$ in $\mathfrak{F}$ after an easy checking of the membership condition for $\Phi'$.

Finally we must show that $\eta$ does not extend to $\eta''$ as stated in the Second Killing Lemma 2.8. Otherwise $\Phi \subseteq \Phi'' : G'' \oplus H'' \longrightarrow R$ and $G'' = x_0 R \oplus G'' \eta''$, hence

$$r = \Phi''(z\eta'', h) \in R.$$

We calculate $r$ differently using continuity of maps:

$$r = \Phi''(z\eta'', h) = \Phi''((\sum_{i \in \omega} q_i r_i w_i)\eta'', h)$$
$$= \Phi''((\sum_{i < n} q_i r_i w_i)\eta'' + q_n(\sum_{i \geq n} q_i (q_n)^{-1} r_i w_i)\eta'', h))$$
$$= \sum_{i < n} q_i r_i \Phi'(w_i \eta, h) + q_n \Phi''((\sum_{i \geq n} q_i q_n^{-1} r_i w_i)\eta'', h)$$

$$\equiv \sum_{i<n} q_i r_i \Phi'(w_i\eta, h) \equiv \sum_{i<n} q_i r_i(w_i\eta\pi_i) \mod q_n.$$

Hence, in the limit $\sum_{n\in\omega} q_n r_n(w_n\eta\pi_n) = r \in R$ which contradicts our choice after Observation 1.1, see (2.10). $\square$

## 3 Construction of the Reflexive modules

Let $\mathfrak{F}^*$ be the set of all those bilinear forms $\Phi$ in $\mathfrak{F}$ which are also finitely covered. Hence $\mathfrak{F}^* \neq \emptyset$ by Corollary 2.5 and $\mathfrak{F}$ in Section 2 can be replaced by $\mathfrak{F}^*$. Suppose that $\Phi_\alpha \in \mathfrak{F}^*$, $\alpha \in \omega_1$ is an ascending, continuous chain of bilinear forms. We will put additional restriction on this choice later on. Hence $\Phi = \bigcup_{\alpha\in\omega_1} \Phi_\alpha$ is a bilinear form on $G \oplus H$ with $G = \bigcup_{\alpha\in\omega_1} G_\alpha$, $H = \bigcup_{\alpha\in\omega_1} H_\alpha$ and $\text{Dom}\,\Phi_\alpha = G_\alpha \oplus H_\alpha$. First we note

**Observation 3.1** $\Phi$ *is non-degenerative.*

**Proof** If $0 \neq g \in G$ then $0 \neq g \in L \subseteq_* G$ for some $L$ of finite rank, hence $L \subseteq_* G_\alpha$ for some $\alpha \in \omega_1$. However $L$ is free and there is $\varphi \in L^*$ with $g\varphi \neq 0$. We also find $h \in H_\alpha$ with $\Phi(\ ,h) \restriction L = \varphi$ from $\Phi_\alpha \in \mathfrak{F}^*$. Hence $\Phi(g,h) = g\varphi \neq 0$. The other case follows by symmetry. $\square$

Now we must recall the definition of the evaluation map $\sigma = \sigma_G : G \longrightarrow G^{**}$ from Section 1. The following maps are obviously related to the evaluation maps.

$$\mathbb{G} : G \longrightarrow H^* \ (g \longrightarrow \Phi(g,\ )) \text{ and } \mathbb{G}' : H \longrightarrow G^* \ (h \longrightarrow \Phi(\ ,h))$$

We claim that

**Lemma 3.2** *The evaluation maps $\sigma = \sigma_G : G \longrightarrow G^{**}$, $\sigma = \sigma_H : H \longrightarrow H^{**}$ and in particular $\mathbb{G} : G \longrightarrow H^*$ and $\mathbb{G}' : H \longrightarrow G^*$ are injective.*

**Proof** If $0 \neq g \in G$, then by Observation 3.1 there is an $h \in H$ such that $\Phi(g,h) \neq 0$. If $\varphi = \Phi(\ ,h) \in G^*$, then $g\varphi = \Phi(g,h) \neq 0$, hence $g\sigma \neq 0$ and $\sigma$ is injective. The other case follows by symmetry. $\square$

The next lemma explains why we want $\mathbb{G}$ and $\mathbb{G}'$ to be surjective.

**Lemma 3.3** *If $\mathbb{G}$ and $\mathbb{G}'$ are surjective, then $\sigma_G$ and $\sigma_H$ are isomorphisms, and hence $G$ and $H$ are reflexive modules.*

**Proof** If $\mathbb{G}$ and $\mathbb{G}'$ are surjective, then $\mathbb{G}$ and $\mathbb{G}'$ are isomorphisms by Lemma 3.2. Hence any $\varphi \in G^{**}$ can be viewed as an element in $H^*$ from $G^* = \mathrm{Im}\,(\mathbb{G}') = \Phi(\ ,H)$ identifying $G^*$ and $H$ under $\mathbb{G}'$. Hence $\varphi \in H^* = \mathrm{Im}\,(\mathbb{G}) = \Phi(G,\ )$ and we find $g \in G$ with

$$\varphi = \Phi(g,\ ). \tag{3.1}$$

Hence
$$\sigma(g)\Phi(\ ,h) = \Phi(g,h) = \Phi(g,\ )(\Phi(\ ,h))$$
for all $h \in H$. Now $\Phi(\ ,h)$ runs through all of $G^*$ and $\sigma(g) = \Phi(g,\ ) = \varphi$ by (3.1). Hence $\sigma_G$ is surjective and an isomorphisms by Lemma 3.2. The proof for $\sigma_H$ is similar, and $G, H$ are reflexive by definition. $\square$

Surjectivity of $\mathbb{G}$ and $\mathbb{G}'$ will be a consequence of the particular choice of the filtrations $\{G_\alpha : \alpha \in \omega_1\}$ and $\{H_\alpha : \alpha \in \omega_1\}$ for $G$ and $H$ respectively. We formulate the restrictions for $H$, the conditions on $G$ follow by symmetry.

If $\varphi \in H_\alpha^*$ is essential for $\Phi_\alpha$, then we can find $\beta > \alpha$ such that

$$\text{any extension } \varphi' \in H_\beta^* \text{ of } \varphi \text{ will not extend to } H_{\beta+1}^*. \tag{3.2}$$

**Lemma 3.4** *(i) If (3.2) holds for $H$, then $\mathbb{G} : G \longrightarrow H^*$ is surjective.*

*(ii) If (3.2) holds for $G$, then $\mathbb{G}' : H \longrightarrow G^*$ is surjective.*

**Proof** It is enough to show that $\mathbb{G}$ is surjective. We will write
$$\mathbb{G}(g) = \Phi(g,\ ) = f_g \in H^* \text{ for all } g \in G.$$

If $\mathbb{G}$ is not surjective, then $\Phi(G,\ ) \neq H^*$ and there is

$$\varphi \in H^* \text{ such that } \varphi \neq f_g \text{ for all } g \in G. \tag{3.3}$$

Obviously we can find $\alpha \in \omega_1$ such that (3.3) 'restricted to $\alpha$' holds, that is

$$\varphi \upharpoonright H_\alpha \in H_\alpha^* \text{ such that } \varphi \upharpoonright H_\alpha \neq f_g \upharpoonright H_\alpha \text{ for all } g \in G_\alpha.$$

Consider $g_1, \ldots, g_n \in G_\alpha$ and $D = \bigcap_{1 \leq j \leq n} \ker f_{g_j} \upharpoonright H_\alpha$, then we can write $H_\alpha = L \oplus D$ with $L$ free of finite rank $\leq n$ and $f_{g_j} \upharpoonright D = 0$ for $j \leq n$. We may assume that the elements $g_1, \ldots, g_n$ are independent, hence $f_{g_1}, \ldots, f_{g_n}$ are independent by the proof of Lemma 3.2. Hence $\operatorname{rk} L = \operatorname{rk} L^* = n$ and we find $a_j \in R$ with $\varphi \upharpoonright L = \sum_{1 \leq j \leq n} f_{g_j} a_j \upharpoonright L$.
If $\varphi \upharpoonright D = 0$, then

$$\varphi = \sum_{1 \leq j \leq n} f_{g_j} a_j = \sum_{1 \leq j \leq n} \Phi(g_j,\ )a_j = \Phi(\sum_{1 \leq j \leq n} g_j a_j,\ ) = \Phi(g,\ ) = f_g$$

from $f_{g_j} \upharpoonright D = 0$ ($j \leq n$) for $g = \sum_{1 \leq j \leq n} g_j a_j$. This contradicts (3.3) 'when restricted to $\alpha$'. Hence any family $g_1, \ldots g_n$ in $G_\alpha$ satisfies

$$\bigcap_{1 \leq j \leq n} \ker \Phi(g_j,\ ) \upharpoonright H_\alpha \nsubseteq \ker \varphi \upharpoonright H_\alpha.$$

There is an $x \in H_\alpha$ such that $x\varphi \neq 0 = xf_{g_i}$ for all $i \leq n$. By an easy rank argument there are many elements $y \in H_\alpha \setminus L$ with $y\varphi = yf_{g_i} = 0$ for all $i \leq n$. So, if $x$ turns out to be in $L$, then replace $x$ by $x + y \in H_\alpha \setminus L$. Hence $\varphi$ is essential for $\Phi_\alpha$. We are able to apply the hypothesis of the First Killing-Lemma 2.7 and see that $\varphi$ does not extend to an elements of $H^*$, hence $\mathbb{G}$ is surjective. □

We have an immediate corollary of Lemma 3.3 and Lemma 3.2

**Corollary 3.5** *If $H$ and $G$ satisfy (3.2) then the following holds for $G$ and $H$.*

(i) $\mathbb{G} : G \longrightarrow H^*$ ($g \longrightarrow \Phi(g,\ )$) and $\mathbb{G}' : H \longrightarrow G^*$ ($h \longrightarrow \Phi(\ ,h)$) are isomorphisms.

(ii) The evaluation maps $\sigma_G : G \longrightarrow G^{**}$ and $\sigma_H : H \longrightarrow H^{**}$ are isomorphisms and $G$ and $H$ are reflexive $R$-modules.

We will need another restriction on the choice of the $\Phi_\alpha$'s taking care of possible isomorphisms between $G$ and $R \oplus G$. This like (3.2) we will do with the help of the prediction principle $\Diamond_{\aleph_1}$, see Eklof, Mekler [5, pp. 139, 140], for instance. We may work in Gödel's universe $V = L$ in which $\Diamond_{\aleph_1}$ holds as shown by R. Jensen. However, there are many other models of set theory with $\Diamond_{\aleph_1}$ in which GCH for instance fails.

Let $X = \bigcup_{\alpha \in \omega_1} X_\alpha$ be an $\aleph_1$-filtration of the set $X$ of cardinality $\aleph_1$ such that at each step $\alpha \in \omega_1$ we have $|X_{\alpha+1} \setminus X_\alpha| = \aleph_0$. If $\Diamond_{\aleph_1}(E)$ holds for some stationary set $E$, then we decompose $E = E_i \cup E_e$ into stationary subsets with $E_i$ in charge of isomorphism and $E_e$ working for essential homomorphisms.

If $\alpha \in E_i$, then let

$$\eta_\alpha : X_\alpha \longrightarrow X_\alpha$$

be the Jensen function predicting maps $X \longrightarrow X$, and if $\alpha \in E_e$, then let

$$\eta_\alpha : X_\alpha \longrightarrow R$$

predict maps $X \longrightarrow R$.

Inductively we define $X_\alpha = G_\alpha \oplus H_\alpha$ and bilinear forms

$$\Phi_\alpha : G_\alpha \oplus H_\alpha \longrightarrow R$$

such that $\Phi_\alpha \in \mathfrak{F}^*$. At limit ordinals $\alpha$ we take unions and it remains to define $\Phi_{\alpha+1}$ for $\Phi_\alpha$. If $\alpha \notin E$ and also if $\alpha \in E_e$ *and* $\eta_\alpha$ is not essential for $\Phi_\alpha$ or $\alpha \in E_i$ *and* $\eta_\alpha \upharpoonright G_\alpha$ is not a monomorphism like $\eta$ discussed in Second Killing Lemma 2.8 and $\eta_\alpha \upharpoonright H_\alpha$ dually for $H$, then we extend $\Phi_\alpha$ trivially to $\Phi_{\alpha+1}$. There are only two interesting cases left.

(i) Suppose $\alpha \in E_e$ and $\eta_\alpha : G_\alpha \longrightarrow R$ is essential for $\Phi_\alpha$ when restricted to $G_\alpha$ or dually $\eta_\alpha : H_\alpha \longrightarrow R$ is essential for $\Phi_\alpha$ when restricted to $H_\alpha$, then we find $\Phi_{\alpha+1}$ from the First Killing-Lemma 2.7 (or its dual version) and kill $\eta_\alpha$. Hence $\eta_\alpha$ does not extend to $H_{\alpha+1}$ or further up.

(ii) Suppose $\alpha \in E_i$ and $\eta_\alpha : G_\alpha \longrightarrow G_\alpha$ is a monomorphism with

$$G_\alpha = x_\alpha R \oplus \operatorname{Im} \eta_\alpha.$$

Then we extend $\Phi_\alpha$ with the aid of Second Killing Lemma 2.8 such that

$$\Phi_{\alpha+1} : G_{\alpha+1} \oplus H_{\alpha+1} \longrightarrow R$$

*and* more importantly $\eta_\alpha$ does not extend to any monomorphism $\eta : G_\beta \longrightarrow G_\beta$ for any $\alpha < \beta$ such that

$$G_\beta = x_\alpha R \oplus \operatorname{Im} \eta.$$

This finishes the construction of $\Phi : G \oplus H \longrightarrow R$, and the following result holds.

**Theorem 3.6** $(ZFC + \diamondsuit_{\aleph_1})$  *There is a reflexive R-module G of cardinality $\aleph_1$ such that $G \not\cong R \oplus G$.*

**Remarks.** Using the Ulam-Solovay's decomposition theorem for stationary sets (see Jech [10]) we may assign different Eklof-invariants $\Gamma(G)$ to various $G$'s, see Eklof, Mekler [5, p. 85] and the example in Theorem 3.6 can be replaced by a family of size $2^{\aleph_1}$ members pairwise non-isomorphic. Also note that the example(s) $G$ are necessarily $\aleph_1$-free as explained earlier.

**Proof of the theorem.** We must check the various consequences from the above construction of $G = \bigcup_{\alpha \in \omega_1} G_\alpha$ *etc.* Clearly the module $G$ is $\aleph_1$-free of cardinality $\aleph_1$. By Corollary 3.5 we must check (3.2) for $G$ and $H$ but this follows from (i) of the construction. Hence $G$ and $H$ are reflexive. If $G \cong xR \oplus G$ then let $\eta : G \longrightarrow G$ be the obvious monomorphism with $G = xR \oplus G\eta$. There is some $\alpha \in \omega_1$ with $x \in G_\alpha$ and $G_\alpha = xR \oplus G_\alpha \eta$ by the modular law. Now we also find a Jensen function $\eta_\beta : G_\beta \longrightarrow G_\beta$ for some $\alpha < \beta$ such that $\eta \restriction G_\beta = \eta_\beta$. At this level $\eta_\beta$ is killed when passing from $\beta$ to $\beta + 1$ by step $(ii)$ of the construction, hence $\eta$ is killed. However this contradicts the fact that $\eta : G_\gamma \longrightarrow G_\gamma$ exists for many $\beta < \gamma$ which follows from a simple back and forth arguments. $\square$

Rüdiger Göbel
Fachbereich 6, Mathematik und Informatik
Universität Essen, 45117 Essen, Germany
e–mail: R.Goebel@Uni-Essen.De
and
Saharon Shelah
Department of Mathematics
Hebrew University, Jerusalem, Israel
and Rutgers University, Newbrunswick, NJ, U.S.A
e-mail: Shelah@math.huji.ae.il